\documentstyle[epsf,stmaryrd,amscd,amstex,amssymb,
11pt
]{amsart}
\input xy

\newcommand{\bu}{\bullet}

\newcommand{\F}{{\cal F}}

\newcommand{\SS}{{\cal S}}

\newcommand{\oplusl}{\bigoplus\limits}

\newcommand{\Gb}{{\overline{G}}}

\newcommand{\Zet}{{\mathbb{Z}}}

\newcommand{\cons}{{\underline {k}}}

\newcommand{\A}{{\sf A}} 

\newcommand{\imbed}{\hookrightarrow}

\newcommand{\la}{\lambda}

\newtheorem{Thm}{Theorem}
\newtheorem{Cor}[Thm]{Corollary}
\newtheorem{Lem}[Thm]{Lemma}
\newtheorem{Prop}[Thm]{Proposition}

\theoremstyle{definition}

\theoremstyle{remark}
\newtheorem{Rem}[Thm]{Remark}
\newtheorem*{Rem*}{Remark}

\numberwithin{Thm}{section}

\newcommand{\proof}{{\em Proof \ }}

\def\proof{\smallskip\noindent {\bf Proof.\ }}

\newcommand\nc{\newcommand}
\nc\on{\operatorname}
\nc\renc{\renewcommand}

\nc{\ppart}{(\!(t)\!)}

\nc\Bun{\on{Bun}}
\renc\mod{\text{-mod}}

\nc\CK{{\mathcal K}}
\nc\CO{{\mathcal O}}

\nc\BC{\mathbb C}
\nc\BA{\mathbb A}
\nc\BO{\mathbb O}

\nc{\fg}{{\mathfrak g}}

\nc{\hg}{\widehat\fg}

\newcommand{\fii}{{\mathfrak i}}
\newcommand{\fj}{{\mathfrak j}}

\numberwithin{equation}{section}

\theoremstyle{definition}

\newcommand{\fC}{{\mathfrak C}}

\newcommand{\Xb}{\bar{X}}

\author{Joseph Bernstein}

\address{\small
School of Mathematical Sciences, Tel Aviv University,
Tel Aviv 69978, Israel
}
\email{bernstei@@math.tau.ac.il}

\author{Roman Bezrukavnikov}

\address{\small
Department of Mathematics, Massachusetts Institute of Technology, 77 Massachusetts ave.,
Cambridge, MA 02139, USA and:
}
\address{\small
National Research University Higher School of Economics,
International Laboratory of Representation
Theory and Mathematical Physics,
20 Myasnitskaya st., Moscow 101000, Russia}
\email{bezrukav@@math.mit.edu}

\author{David Kazhdan}

\address{\small
Einstein Institute of Mathematics,
The Hebrew University of Jerusalem,
Givat Ram, Jerusalem, 9190401, Israel}
\email{kazhdan@@math.huji.ac.il}

\title[Deligne-Lusztig duality and compactification]{Deligne-Lusztig duality and wonderful compactification}

\begin{document}


\begin{abstract}
We use geometry of the wonderful compactification to obtain a new proof of the relation between Deligne-Lusztig (or Alvis-Curtis) duality for $p$-adic groups and homological duality.
This provides a new way to introduce an involution on the set of irreducible representations of the group
which has been defined by A.~Zelevinsky for $G=GL(n)$ and by A.-M.~Aubert  in general
(less direct geometric approaches to this duality have been developed earlier 
by Schneider-Stuhler and by the second author).
As a byproduct, we describe the Serre functor for representations of a $p$-adic group.
\end{abstract}

\maketitle

\centerline{\em{To Sasha Beilinson, with admiration and best wishes for his birthday.}}

\section{Introduction}

The goal of this note is to present a new conceptual proof of the relation between homological
and Deligne-Lusztig dualities and apply it to define an involutive auto-equivalence of the category of admissible representations
which induces an involution on the set of irreducible representations.

Let $G$ be a reductive $p$-adic group. For every smooth $G$-module $M$ one can form
a complex 
$$0\to M \to \oplusl_P 
i_P^Gr_P^G(M)\to \cdots \to i_B^Gr_B^G(M)\to 0,$$
where  $i_P^G$, $r^G_P$ denote, respectively, the parabolic induction and Jacquet 
functors,\footnote{\label{f1} For uniformity of notation we let $i$ and $r$ denote {\em normalized}
or unitary parabolic induction and Jacquet functors. Note that replacing both $i$ and $r$ by ordinary
(rather than normalized) induction and Jacquet functors does not affect the composition
$i_P^Gr_P^G$.}
and summation in the $i$-th term runs over conjugacy classes of parabolic subgroups of corank $i$.
The differential is the sum of the natural maps coming from adjunction between induction and Jacquet functors and transitivity of parabolic induction, taken with appropriate signs. We call this complex
the Deligne-Lusztig complex associated to $M$, and denote it by $DL(M)$. Analogous complexes for  representations of a finite Chevalley group have been considered in \cite{DL}, the corresponding
automorphisms of the Grothendieck group has been studied earlier in \cite{Al}, \cite{Cur}
for finite Chevalley groups and in \cite{Ka} for $p$-adic groups.

Starting from a complex of smooth representations $M^\bu$ one gets a bicomplex $DL(M^\bu)$, thus
we get a functor $DL:D^b(Sm)\to D^b(Sm)$.

The main result of this  note is as follows.

\begin{Thm}\label{main}
For a complex $M$  with admissible cohomology  we have a canonical quasi-isomorphism
$$DL(\check{M})\cong RHom_H(M,H)[r],$$
where $\check{M}$ denotes the contragradient representation, 
$H\cong C_c^\infty(G)$ is the regular
bimodule for the Hecke algebra and $r$ is the split rank of the center of $G$.
\end{Thm}

 Theorem \ref{main} is proved in section \ref{geom}
by using an explicit resolution for the regular bimodule over $H$ coming from geometry of 
the wonderful compactification $\Gb$ of $G$. The idea of this proof is as follows.

Recall that $\Gb$ admits a stratification where the strata are indexed by conjugacy classes of parabolic subgroups in $G$. A sheaf on a stratified space admits a standard resolution with terms indexed by the strata. We apply this to the sheaf of smooth sections (see \ref{nota} for the definiton)
of the sheaf $j_*(\cons_G)$, where $j:G\to \Gb$ is the embedding and $\cons_G$ is the constant sheaf on $G$. Taking global sections, we get a resolution for the regular bimodule over $H$.

Next, specialization to the normal bundle for functions on a $p$-adic manifold \cite{BK} allows us to describe the terms of this resolution via functions
on the normal bundles to the strata. A normal bundle to a stratum $\Sigma$ contains an open $G\times G$ orbit $O_\Sigma$, and functions on $O_\Sigma$ 
form an injective $G$-module, dual to the one appearing in the study of second adjointness \cite{BK}. Taking $Hom$ from a finitely generated module $M$ into the space of functions on $O_\Sigma$
yields a term of the complex $DL(\check{M})$. 
On the other hand, a term of the resolution for $H$ described in the previous paragraph receives a surjective (but not injective) map from the space of functions on $O_\Sigma$.
 It turns out that this map induces an isomorphism on $Ext$'s from an 
  admissible module. This yields Theorem \ref{main}.

 Other less direct geometric proof of Theorem \ref{main} based on Borel-Serre compactification of the Bruhat-Tits building of $G$
 have appeared in \cite{SS} and  the PhD thesis of the second author \cite{Bez}.
 
In section \ref{alg} we present a generalization and algebraic consequences of Theorem \ref{main}. 
This section (except for \S \ref{Serre_f} which describes the Serre
functor for admissible modules over a $p$-adic group) follows the strategy devised by the first author 
several decades ago,  see announcement in \cite{Ber2nd}.
  
The main application recorded in section \ref{alg}
is Corollary \ref{CorAub} which says that for an irreducible module $M$
the complex $DL(M)$ has cohomology in only one degree $d$, where $d$ is dimension of the component of Bernstein center containing $M$.  This implies that in fact $H^d(DL(M))$ is an irreducible
representation, thus we get an automorphism on the set of representations. This result has already 
appeared in \cite{Aub} and in \cite{SS}. 
In fact, the automorphism on the set of irreducible representations can be shown to be an involution, which for $G=GL(n)$ has been defined and studied in \cite{Ze},
thus it is often called the Aubert, or the Zelevinsky involution.
While the proof in \cite{Aub} is short, Aubert-Zelevinsky involution has many interesting known or expected properties related to Langlands conjectures,\footnote{Namely, 
it is expected to interchange the Arthur and the Deligne-Langlands parameters, see \cite{Hi} for details.} Koszul duality (cf. \cite{MR}) etc., which are far from being completely 
understood. The central idea of the present note is a link between this involution and geometry of the wonderful compactification, which is, on the other hand, related, as we plan to argue elsewhere, to the local trace formula; see also \cite{DW}, 
 \cite{Ga} for other apparently related constructions. We hope this link will help clarify some of the outstanding questions about the involution.
 
Finally, let us remark that a $p$-adic group $G$ of rank 1 is known to be a hyperbolic group in the sense of Gromov \cite{Gr}, the wonderful
compactification can be shown to be closely related to the Gromov compactification $\bar{G}_{hyp}$ of $G$ viewed
as a hyperbolic group (namely, $\bar{G}$ is the closure of the image of the antidiagonal embedding of
$G$ into $\bar{G}_{hyp}^2$). It is tempting to look for a generalization of the results of this note to other
hyperbolic and semi-hyperbolic (in the appropriate sense) groups.
 
 {\bf Acknowledgements.}  
 We thank Vladimir Drinfeld for many helpful conversations over the years.
 The second author is also grateful to  Michael Finkelberg, Leonid Rybnikov and Jonathan Wang
 for  motivating discussions. The impetus for writing this note came from a talk given by the second author at the Higher School for Economics (Moscow), he thanks this institution for the stimulating opportunity.  Finally, we thank Victor Ginzburg for an inspiring correspondence which has motivated
 section \ref{Serre_f}.

The project  received funding from ERC under grant agreement No 669655.
J.B. acknowledges partial  support by the ERC grant 291612,
R.B. was partly supported by the NSF grant DMS-1601953
and   Russian Academic Excellence Project '5-100,  The collaboration of R.B. and D.K. was supported by the US-Israel BSF grant 2016363.

\section{Proof of Theorem \ref{main}} \label{geom}

\subsection{Notation}\label{nota} Let $k$ denote an algebraically closed characteristic zero field. 

Let also $F$ be a local non-Archimedean field and $G$ be the group of $F$ points of a 
reductive
algebraic group of $F$. We let $\Gb$ denote its wonderful, or De Concini -- Procesi (partial) compactification.
The definition of wonderful compactification is standard in the case of an adjoint semi-simple group only. For a general semi-simple group $\Gb$ is defined to be the normalization of $\Gb_{ad}$ in $G$;
for a general reductive group we let $\Gb= (\Gb' \times C)/(G'\cap C)$, where $C$ is the center of $G$,
$G'$ is the derived group and $G_{ad}=G/C$. Thus $\Gb$ is proper if $G$ is semi-simple,
while $\Gb(F)$ is compact if $G$ has compact center. From now on we will abuse the
notation by using the same notation for an algebraic variety over $F$ and the 
topological space of its $F$ points.

For a totally disconnected topological space $X$ we write  $C_c(X)$ for the space of compactly supported functions on $X$, $C^\infty(X)$ for the space of locally constant function and   $C_c^\infty=C_c(X)\cap C^\infty(X)$; all functions are assumed to be $k$-valued.

We have the {\em Hecke algebra} $H$  which is the space of locally constant compactly supported measures on $G$ equipped with the convolution algebra structure; a choice of a Haar measure yields an isomorphism $H\cong C_c^\infty(G)$. For an open compact subgroup $K\subset G$ we have
the unital subalgebra $H(G,K)\subset H$ of $K$-biinvariant measures.

Recall that for a module $M$ over a totally disconnected group a vector $v\in M$ is called smooth if
its stabilizer is open. 
If $M$ is a module over an open subgroup in $G$ we let $M^{sm}$ denote the submodule of smooth vectors. If $M=M^{sm}$ then $M$ is called smooth.

We let $Sm(G)$ denote the category of 
smooth $G$-modules. For $M\in Sm(G)$ the  
 contragradient module $(M^*)^{sm}$ is denoted by $\check{M}$ or $M\check{\ }$.
 We let $Adm(G)\subset Sm(G)$ be the subcategory of admissible representations.

A $G$-equivariant sheaf $\F$ on a $G$-space $X$ will be called smooth if for every point $x\in X$ and every local section $s$ of $\F$ defined on a neighborhood of $x$ there exists an open subgroup
$K$ in $G$ and a $K$-invariant neighborhood $U$ of $x$, such that $s|_U$ is (well defined and)
$K$-invariant. For a $G$-equivariant sheaf $\F$ we 
let $\F^{sm}$ denote the maximal smooth subsheaf. In other words, $s\in \Gamma(U,\F)$ belongs to 
$\F^{sm}$ if $U$ admits an open covering $U_i$ such that $U_i$ is invariant under an open subgroup $K_i$ and $s|_{U_i}$ is $K_i$-invariant.

It is easy to see that $\Gamma(X,\F)$ is a smooth $G$-module provided that $\F$ is a smooth $G$-equivariant sheaf on a compact $G$-space $X$.

The constant sheaf with stalk $M$ over a topological space $X$ will be denoted by $\underline{M}_X$.

\subsection{The specialization complex}
We start by introducing an auxiliary complex of $G\times G$ modules.

For a subset $S$ in the set  $I$ of simple roots let $X_S=(G/U_S\times G/U_S^-)/L_S$, where 
$P_S=L_S U_S$ is a parabolic subgroup in the conjugacy class determined by $S$ and $P_S^-=L_S U_S^-$ is an opposite parabolic. In \cite{BK} we have defined\footnote{In \cite{BK}
we worked under the running assumption that the group $G$ is split. This assumption is not used in proving any of the facts referenced in the present paper.}
 the cospecialization map
$c_S: C_c^\infty(X_S)\to C_c^\infty (G)$ for all $S$. Similar considerations yield a map 
$c_{S,S'}:C_c^\infty(X_S)\to  C_c^\infty(X_{S'})$ for $S\subset S'$, such that $c_{S',S''}c_{S,S'}=c_{S,S''}$ and $c_{S,I}=c_S$ which we also call the cospecialization map.

Taking the sum of these maps with appropriate signs we get a complex of $G\times G$ modules:
$$0\to  C_c^\infty(X_I) \to \cdots \to  \oplusl_{|S|=1}  C_c^\infty(X_S)\to  C_c^\infty(G) \to 0,$$
%
which we call the cospecialization complex. Passing to the contragradient $G\times G$-modules we get
a complex which we call the specialization complex $\fC$.

\begin{Lem}
a) The terms of the specialization complex $\fC$ are injective as left $G$-modules.

b) For a finitely generated $G$-module $M$ the complex
 $Hom(M,\fC)$ is canonically isomorphic to  the complex $DL(\check{M})$.
\end{Lem}

\proof It is clear that 
$$C_c^\infty(X)=i_{P\times P^-}(C_c^\infty(L)) \cong i_P^G(C_c^\infty(G/U^-)),$$
 where $i$ denotes unitary parabolic induction.\footnote{Notice that $X_S$ has a $G$-biinvariant measure,
 while $G/U$ has a $G\times L$ invariant measure, thus replacing the unitary parabolic induction here by
 the ordinary does not affect the result (up to an isomorphism), cf. footnote \ref{f1}.}
Since passing to the contragradient module commutes with unitary parabolic induction, we see that the terms of $\fC$ are
of the form 
$$i_{P\times P^-}^{G\times G}((C_c^\infty(L)^*)^{sm}) = i_{P}^{G}(C(G/U^-)^{sm}),$$ 
where in the last expression the superscript denotes taking smooth vectors with respect to 
the $G\times L$ action, and compactness of $G/P^-$ is used. 

For an open compact subgroup $K\subset G$ the space $(C(G/U^-)^{sm})^K$ is 
a direct sum over $K$-orbits on $G/P^-$ of $L$-modules of the
form $C(K_L \backslash L)$ for an open compact subgroup $K_L \subset L$.
Since $C(K_L \backslash L)$ is an injective object in $Sm(L)$, it follows that
$C(G/U^-)^{sm}$ is also an injective 
object in $Sm(L)$.
 The standard Frobenius adjointness shows that parabolic induction sends injective modules to injective ones,  this proves (a). 

For  $M\in Sm(G)$ we have:
$$Hom_G(M, i_{P}^{G}(C(G/U^-)^{sm}))= Hom_L(r^G_P(M), C(G/U^-)^{sm} ) .$$
If $M$ is finitely generated then $r^G_P(M)$ is also finitely generated, and for a finitely generated
smooth $L$-module $N$ we have 
$$Hom_L(N, C(G/U^-)^{sm} )\cong i_{P^-}^G (Hom (N, C(L)^{sm}))=i_{P^-}^G(\check{N}).$$
Recalling the isomorphism $r^G_P(M)\check{\ }\cong r^G_{P^-}(\check{M})$ which follows directly from  the second adjointness, we get statement (b).   \qed

\medskip

We do not describe cohomology of $\fC$; instead, in the next subsection we define a certain quotient
of $\fC$ which will be shown to be a resolution for $C_c^\infty(G)$.

\subsection{The resolution}

For a subset $S$ in the set  $I$ of simple roots let $\Gb_S$ be the corresponding closed stratum in 
 $\Gb$ and let $i_S:\Gb_S\to \Gb$, $j:G\to \Gb$ be the embeddings.

 We let $\F_S=i_S^* (j_*(\cons_G)^{sm})$. 

We set  $M_S:=\Gamma(\F_S)$. In particular, $M_\emptyset$ is the space of sections of the 
constant sheaf on $G$, i.e. the space $C^\infty(G)$ of locally constant functions on $G$.

We now proceed to construct a resolution of the regular bimodule for $H$. This will be derived
from a standard resolution for a sheaf on a stratified space.

Let $X$ be a 
topological space and $Z_1,\dots, Z_n$ be closed subspaces.
For a subset $S\subset [1,n]$ let $X_S$ denote the intersection of $Z_i$, $i\in S$ (in particular, $X_\emptyset=X$) and $i_S:X_S\to X$ be the embedding; let $j:U\to X$ be the embedding of 
the complement to $\cup_i Z_i$. Then for a sheaf $\F$ on $X$ we can form a
complex of sheaves on $X$:

\begin{equation}\label{vkl_iskl}
0\to (i_\emptyset)_*i_\emptyset^*(\F) \to \cdots \to \oplusl_{|S|=i} (i_S)_*i_S^*(\F)\to \cdots
\to (i_{[1,n]})_*i_{[1,n]}^*\F \to 0,
\end{equation}
where the differential is the sum of restriction maps with appropriate signs.

\begin{Lem}\label{Lem_vkl_iskl}
The complex \eqref{vkl_iskl} is exact in positive degree and its cohomology in degree zero is  $j_!j^*(\F)$.
\end{Lem}

\proof A complex of sheaves is exact if an only if the complex of stalks at every
point is exact. Applying this criterion to the complex  
$$ 0\to j_!j^*(\F)\to (i_\emptyset)_*i_\emptyset^*(\F) \to \cdots \to \oplusl_{|S|=i} (i_S)_*i_S^*(\F)\to \cdots \to (i_{[1,n]})_*i_{[1,n]}^*\F \to 0$$
we get the claim. \qed

\begin{Cor}\label{resol}
We have  a complex of $G\times G$ modules
\begin{equation}\label{compl}
0\to M_\emptyset \to \oplusl_{|S|=1} M_S \to\cdots \to M_I\to 0,
\end{equation}
whose cohomology in concentrated in degree zero.

If $G$ is semisimple, then the 0-th cohomology of this complex is isomorphic to $C_c^\infty (G)$.\end{Cor}

\proof  We apply  Lemma \ref{Lem_vkl_iskl} to  the sheaf $j_*(\underline{k})^{sm}$ on $\Gb$. 
Since $\Gb$ is totally disconnected, the sheaves have no higher cohomology, so the resolution for a sheaf yields a resolution for its global sections. It is clear that if $\Gb$ is compact, then $C_c^\infty(G)=\Gamma_c(\underline{k}_G) = \Gamma(j_!(\underline{k}_G))$.
\qed

Our next goal is  a more explicit description of the complex \eqref{compl}.

If $G$ is adjoint, then $\Gb$ is  (the space of $F$-points of) a nonsingular projective algebraic variety.
Let $\Xb_S$ denote the normal bundle to $\Gb_S$ in $\Gb$. For a general $G$ we let $\Xb_S$
be the quasi-normal cone in the sense of \cite[\S 11.1]{BK}.
Let $\fii_S$ be the zero section embedding of $\Gb_S$ into $\Xb_S$.

Recall that  $\Xb_S$ contains an open dense $G\times G$-orbit $X_S$ isomorphic
to $(G/U_S\times G/U_S^-)/L$. The projection $X_S\to G/P\times G/P^-$ extends to a map
$\Xb_S\to G/P\times G/P^-$; let $L^+$ denote a fiber of that map. Thus 
 $L^+$ is partial compactification of $L$.
 
 For future reference we recall the following description of $L_+$. Let $C$ be the center of $L$
 and $C_+$ be its partial compactification defined as $Spec(k[\Lambda_C^+])$ where $\Lambda_C$
 is the character lattice of $C$ and $\Lambda_C^+\subset \Lambda_C$ is the semigroup of dominant weights. Then we have
 \begin{equation}\label{Lp}
 L_+=(\overline{L}\times C_+)/C,
 \end{equation}
 where $\overline{L}$ is the (partial) wonderful compactification of $L$.

 Let $\fj_S:X_S\to \Xb_S$ denote the embedding.

\begin{Prop}\label{FSMS}
a) We have  canonical $G\times G$-equivariant isomorphisms:
 \begin{equation}\label{FS}
 \F_S\cong \fii_S^*(\fj_{S*} (\cons_{X_S
 })^{sm});
 \end{equation}
 \begin{equation}\label{MS}
M_S=i_{P_S\times P_S^-}^{G\times G} (( C (L)^{sm}/C_0 (L)^{sm}),
\end{equation}
  where $C_0(L)$
is the space of  
functions on $L$ whose support is closed in $L^+$. 

b) The surjections $i_{P_S\times P_S^-}^{G\times G} (C (L)^{sm})=C(X_S)^{sm}\to M_S$
from \eqref{MS} provide a morphism of complexes of $G\times G$-modules from the specialization
complex $\fC$ to the complex \eqref{compl}. 
\end{Prop}

\proof  We first show that isomorphism \eqref{MS} follows from \eqref{FS}.  
We have an exact sequence of sheaves 
on $X_S$:
$$0\to \fj_{S!} (\cons_{X_S
})^{sm}\to
\fj_{S*} (\cons_{X_S
})^{sm}\to \fii_S^*(\fj_{S*} (\cons_{X_S
})^{sm})\to 0.$$
It is clear that
 $$\Gamma( (\fj_{S*}\cons_{X_S
 })^{sm}) \supset  \Gamma(\fj_{S*}\cons_{X_S
 }) ^{sm}=C(X_S
 )^{sm}= i_{P_S\times P_S^-}^{G\times G} (C(L)^{sm}),$$
  $$ \Gamma(\fj_{S!} (\cons_{X_S
  }) )^{sm} = i_{P_S\times P_S^-}^{G\times G} (C_0 (L)^{sm}).$$
 It is also easy to see that the functor $\F \mapsto \Gamma(\F)^{sm}$ is exact on the category of smooth $G$-equivariant sheaves on a totally disconnected $G$-space, this yields \eqref{MS}.

We proceed to check \eqref{FS}. We deduce it from the next general Lemma, to state it we need some notation.

Let $X$ be a smooth algebraic variety over $F$ and $D\subset X$ be a divisor with normal crossings.
Let $D_1,\dots, D_n$ be the components of $D$; for $I\subset [1,n]$ let 
$D_I=\bigcap\limits_{s\in I} D_s$, let $i_I:D_I\imbed X$ be the embedding and let $D_I^o\subset D_I$ be the complement to the union of $D_J$, $J\supsetneq I$. Set $U=D_\emptyset ^o$ and let $j:U\to X$ be the embedding.

Let $N_I=N_X(D_I)$ be the normal bundle to $D_I$ in $X$ and $\fii_I:D_I\to N_I$ be the zero section embedding
and let $\fj_I:U_I\imbed N_I$ be the embedding, where $U_I$ is the complement to the union
of $N_{D_i}(D_I)$, $i\in [1,n]$. 

\begin{Lem}\label{25}
Suppose that an algebraic group $G$ acts on $X$ preserving $D_i$, so that the action on $D_I^o$
is transitive for each $I$.  Assume\footnote{This assumption may in fact be redundant.} also that for every $I$ and for any (equivalently some) point $x\in D_I^o$ the action of the stabilizer $Stab_G(x)$ on the fiber at $x$ of the normal bundle to $D_I$ in $X$ decomposes as a sum of linearly independent 
characters.\footnote{It is easy to see that the running assumptions imply it
decomposes as a sum of one-dimensional representations.}

Then for each $I$ we have a canonical $G$-equivariant isomorphism 
of sheaves on $D_I$:
$$i_I^* j_*(\cons_U)^{sm}\cong \fii^*\fj_{I*}(\cons_{U_I})^{sm}.$$  
\end{Lem}

\proof 
Choose a neighborhood $V$ of $D_I$ in $X$ and a map $\phi$
from $V$ to $N_I$ which is admissible in the sense of
\cite[Definition 3.2]{BK}. It is clear that $\phi$ induces an isomorphism 
$\iota_\phi: i_I^*j_*(\cons_U)\to \fii^*_I \fj_*(\cons_{U_I})$.
It is also easy to see that $\iota_\phi$ restricts to an isomorphism 
$\iota_\phi^{sm}: i_S^*j_*(\underline{k})^{sm}\to \fii^*_S \fj_{I*}(\underline{k})^{sm}$. 

We claim that $\iota_\phi^{sm}$ (though possibly not $\iota_\phi$) is independent of the choice of $\phi$. We deduce this from  \cite[Lemma 3.4]{BK}. 

Let $\phi'$ be another such map. 
Notice that the bundle $N_I$ splits as a sum of line bundles $N_{D_J}(D_I)$, where $J\subset I$ and $|J|=|I|-1$. Thus the torus $(F^*)^d$, $d=|I|$ acts
on $N_I$. Fixing a uniformizer $\varpi\in F$ we get an embedding $\Zet^d\imbed (F^*)^d$,
$\la\mapsto \varpi^\la$.
According to {\em loc. cit.},
for any locally constant function $f$ with compact support on  $U_I$ 
 there exists $\la_0\in \Zet^d$, such that $$\phi_*(\varpi^\la_*(f))=\phi'_*(\varpi^\la_*(f))$$
for $\la\geq \la_0$, where "$\geq$" refers to the coordinatewise partial order on $\Zet^d$, i.e. $\la=(\la_1,\dots, \la_d) \geq \mu = (\mu_1,\dots, \mu_d)$
if $\la_i\geq \mu_i$ for all $i=1,..,d$. 

It is clear that for $f$ as above and a Taylor series in $d$ variables $\alpha=\sum a_{\la}t^{\la}$
the function $\alpha(f):=\sum a_\la \varpi^\la_*(f)$ is well defined, the above implies that 
$$\phi_*(\alpha(f))=\phi'_*(\alpha(f))$$
for $\alpha\in t^{\la_0}k[[t]]$. 

Let now $K$ be an open compact subgroup in $G$ and $C$ a compact open $K$-invariant subset in $D_I^o$. The second assumption of Lemma implies that $G$ acts transitively on $U_I$.
It is easy to deduce that there exist a finite collection of functions $f_i$ with compact support 
on the preimage of $C$ in $U_I$, 
 such that the sections of the form $\fii_I^*(\alpha(f_i))$, $\alpha\in k[[t]]$ span the space of $K$-invariant sections of $\fii_I^*(\fj_I)_*(\cons_{U_I
 })|_C$. Furthermore, this remains true if we only allow $\alpha\in t^{\la_0}k[[t]]$ for a fixed $\la_0\in \Zet^d$. Now we see that
$\phi_*(\sigma)=\phi'_*(\sigma)$ for any $K$-invariant section of $\fii_I^*(\fj_I)_*(\cons_{U_I
})|_C$,
thus we have checked independence of $\iota_\phi^{sm}|_{D_I^o}$ on $\phi$. 

We finish the proof by checking
that the  isomorphism on stalks at every point in $x\in D_I$ 
induced by $\iota_\phi^{sm}$ is independent of $\phi$. Let $x$ be such a point.
The statement is equivalent to saying that an automorphism of a neighborhood $V$
of $D_I$ in $X$ which is identity on $D_I$, preserves $D_J\cap V$ for all $J$ and has normal component of differential at points of $D_I$ equal identity acts by identity on the stalk of 
$i_I^*j_*(\cons_U)$ at $x$. In view of the previous paragraph we already know
this is true for points $x\in D_I^o$. Let now $x\in D_I$ be arbitrary, we have $x\in D_J^o$ for some $J$. Applying the argument
of the previous paragraph to $D_J$ and observing that an automorphism of a neighborhood of $D_I$ 
 as above restricts to an automorphism of a neighborhood of $D_J$ ($J\supset I$) satisfying similar conditions, we get independence of $\iota_\phi^{sm}$ on $\phi$.

 Independence of $\phi$ also shows that the isomorphism is 
compatible with the natural $G$-equivariant structures on the two sheaves.
This proves the Lemma. \qed

We return to the proof of   Proposition \ref{FSMS}. Isomorphism \eqref{FS} in the adjoint case follows directly from Lemma \ref{25}.
In the general case it is constructed in a similar way, relying on a generalization of Lemma \ref{25}
to the setting of \cite[Appendix]{BK}, the proof of this generalization is parallel to the proof of Lemma \ref{25}. This proves part (a) of the Proposition. 

 Compatibility with the differential
claimed in part (b) follows by comparing definition of $c_{S,S'}$ from \cite{BK} with the proof of the Lemma.  \qed

\begin{Prop}\label{ext_van}
For an admissible $L$-module $M$ we have
$$Ext^{i}_L(M, C_0^\infty (L))=0 $$
for all $i$.
\end{Prop}

\proof  
Let $g\in C$ be an element which has a positive pairing with elements of $C_+$. We claim that for a nonzero constant $c$ the element $(c-g)$ acts on $C_0^\infty(L)$ by an invertible
operator, this yields the Proposition since for an admissible module $M$ Schur Lemma shows that
$\prod (c_i-g)$ kills $M$ for some nonzero constants $c_1,\dots,c_n$.

The operator inverse to $c-g$ is obtained by expanding $(1-c \cdot g^{-1})^{-1}$ in a series
$1+cg^{-1}+c^2g^{-2}+\dots$. In order to check that this expression for the inverse operator is well
defined it suffices to check that for a subset $X$ in $L$ which is closed in $L_+$ and for a point
$x\in L$ there exist only finitely many natural numbers $n$ for which $g^{-n}(X)\owns x$. 
This is clear from \eqref{Lp}
which implies that the sequence $g^n(x)$ has a limit in $\partial L=L_+\setminus L$.
 \qed

\subsection{The main result}
We are now ready to prove:

\begin{Thm} \label{main1}
Let $r$ be the split rank of the center of $G$.
%
For a complex $M$ with admissible cohomology we have a canonical isomorphism in the 
derived category:
$$RHom(M, C_c^\infty(G))\cong DL(\check{M})[-r].$$
\end{Thm}

\proof It is easy to reduce the statement to the case of a group $G$ with compact center.
In that case it follows by
comparing Corollary \ref{resol} and Proposition \ref{FSMS}(b) with Proposition \ref{ext_van}. \qed

\section{A generalization and algebraic consequences} \label{alg}

\subsection{Cohen-Macaulay property}
Let $\A$ be an indecomposable summand in the category $Sm(G)$. According to \cite{cen},
the center $Z$ of $\A$ (i.e. the endomorphism ring of the identity functor $Id_\A$) is the ring
of functions on an algebraic torus invariant under an action of a finite group; in particular, $Z$
is a {\em Cohen-Macaulay} commutative algebra. 

By \cite{cen} 
the category $\A$ is equivalent to the category of 
$H_\A$-modules for a certain $Z$-algebra $H_\A$, which is finitely generated as a $Z$-module. The algebra $H_\A$ is not uniquely defined, but it is defined uniquely up to a (canonical) Morita equivalence. 

Notice that given a Morita equivalence between two $k$-algebras $R$ and $R'$ one also gets a Morita 
equivalence between $R\otimes R^{op}$ and $R'\otimes (R')^{op}$ sending the regular bimodule
for $R$ to the regular bimodule for $R'$. This shows that 
 the homological duality functor $M\mapsto RHom_R(M,R)$ considered as a (contravariant) functor from 
 the derived category of $R$-modules to the derived category of $R^{op}$-modules commutes with Morita equivalences; it is clear that
for $R=H_\A$ we recover the above homological duality.

\begin{Prop}\label{H_CM}
$H_\A$ is a Cohen-Macaulay $Z$-module.
\end{Prop}

\proof It suffices to show that a projective generator for $\A$ is a Cohen-Macaulay $Z$-module.
In view of second adjointness \cite{Ber2nd}, \cite{BK}, a parabolic induction functor sends projective objects to projective ones. Recall that to the indecomposable summand $\A$ there corresponds a Levi
subgroup $L$ in $G$ and a cuspidal representation $\rho$ of $L$, which is defined uniquely up to twisting with an unramified character of $L$. Let $L^0\subset L$ be the kernel of unramified characters, thus $L^0$ is the subgroup generated by all compact subgroups. Let $P=L\cdot U$ be a parabolic. Then the module $\Pi_\rho =  i_P^G(ind_{L^0}^L (\rho|_{L^0}))$ is a projective generator for
 $\A$. The action of $Z$ on $\Pi_\rho$ extends to an action of the ring $\tilde Z = k[L/L^0]$, the ring
 of regular functions on the algebraic torus. Moreover, $\Pi_\rho$ is easily seen to be a free module over $\tilde Z$, hence it is a Cohen-Macaulay module over $Z$. \qed

Following \cite{BBG}, by a (possibly noncommutative) {\em Cohen-Macaulay} ring of dimension $d$ we will understand
a ring $A$ for which there exists a homomorphism from a commutative ring $C$ to the center of $A$,
such that $C$ is a Cohen-Macaulay ring  of dimension $d$ and $A$ is a finitely generated Cohen-Macaulay module over $C$ of full dimension. Standard properties of Cohen-Macauliness show that when $A$ is commutative this agrees with the standard definition.

The following statement is standard for commutative rings, the noncommutative generalization is proved in a similar way.

\begin{Lem}\label{Extd}  Let $A$ be a Cohen-Macaulay (possibly non-commutative) ring of dimension $d$ and $M$ a finite dimensional  
module. Then $Ext_A^i(M,H)=0$ for $i<d$.
\end{Lem}

\proof We use induction in $d$. Let $C\to A$ be as above. We can find a non-zero divisor $f\in C$ which
annihilates $M$. Then $\bar{A}=A/fA$ is a Cohen-Macaulay ring of dimension $d-1$ and 
$$Ext^i_A(M,A)=Ext^i_{\bar{A}}(M, RHom_A(\bar{A},A))=Ext^{i-1}_{\bar{A}}(M,\bar{A}),$$
where we used that $Ext^i_A(\bar{A},A)=\bar{A}$ for $i=1$ and this Ext vanishes for $i\ne 1$. \qed

\subsection{Aubert-Zelevinsky involution}

\begin{Thm}\label{CorAub} 
Let $r$ denote the split rank of the center of $G$.

a) Let $M$ be an admissible module belonging to a component of Bernstein center of dimension $d$.
Then $DL(M)$ has cohomology in degree $d-r$ only.

b) The functor $M\mapsto H^{d-r}(DL(M))$ is an autoequivalence of the summand in $Adm(G)$ corresponding to a component in the Bernstein center of dimension $d$.
\end{Thm}

\proof a) 
 The complex $DL(M)$ is concentrated in degrees from 0 to $d-r$. On the other hand,
Proposition \ref{H_CM} and Lemma \ref{Extd} show 
 that $Ext_H^i(\check{M}^K,H_\A)=0$ for $i\leq d-r$.

b) The functor $DL$ is an autoequivalence of the derived category
 $D^b(Sm(G))_{adm}$ (complexes with admissible cohomology). Part (a) shows that on a summand
 of dimension $d$ this autoequivalence shifted by $d-r$ preserves the abelian heart.

\begin{Cor}\label{DLCor}
 Given an irreducible module $M$ belonging to a component of dimension $d$  there exists an irreducible module $M'$ such that:
$[DL(M)]=(-1)^{d-r} [M']$.
\end{Cor}

\begin{Rem}
See also \cite[III.3 and Proposition IV.5.1]{SS} for statements equivalent to Theorem \ref{CorAub}
and Corollary \ref{DLCor}.
\end{Rem}

\subsection{Homological duality and Deligne-Lusztig duality for non-admissible modules}
Notice that both the homological duality and $DL$ are defined on the  derived category
of finitely generated smooth modules, so it is natural to ask if they can be related as functors on that larger category.

To state the answer we introduce the {\em Grothendieck-Serre} duality $\SS$ on $Sm$.

If $A$ is noncommutative ring over a commutative ring $R$, then for an $A$-module $M$ and 
a $R$-module $N$ the space $Hom_R(M,N)$ carries a right $A$-module structure.
Passing to derived functors we get a functor $RHom:D^b(A-mod)^{op}\times D^+(R-mod)\to D^+(A^{op}-mod)$.

We now consider the category $Sm_R$ of smooth $G$-modules over a base commutative ring $R$
which is finitely generated over $k$.

Since we have a standard isomorphsim $H^{op}\cong H$, we get a functor $RHom:D^b(Sm_R)^{op}\times D^+(R-mod)\to D^+(Sm_R)$.

\begin{Prop}\label{RHom_N}
Let $D^b(Sm_R)_{adm}$ be the full subcategory in $D^b(Sm_R)$ consisting of complexes whose cohomology is $R$-admissible.
Then for $M\in D^b(Sm_R)_{adm}$ and a $N\in D^b(R-mod)$ we have a natural
isomorphism $$RHom_{D^b(Sm_R)}(M,H\otimes N)\cong DL(RHom_R(M,N)).$$ 
\end{Prop}

\proof The proof is similar to the proof of Theorem \ref{main}.
Let $\fC_N$ denote the specialization complex tensored by $N$ over $k$, this is a resolution for
$H\otimes N$. It is clear that
$$RHom _{Sm_R} (M, \fC_N)\cong DL(RHom_R(M,N)).$$
Thus it remains to generalize Proposition \ref{ext_van} by checking that for a parabolic $P=LU$ we have:
$Ext^i_{Sm_R(L)}(r^G_P(M),  C_0^\infty (L)\otimes N)=0 $ for all $i$. 
Let $g\in C$ be as in the proof of Proposition \ref{ext_van}.
Since 
$r^G_P(M)$ is an $R$-admissible finitely generated $L$-module, 
 we can find a monic polynomial $P\in R[t]$, such that  $P(g)$ annihilates $r^G_P(M)$.
 It remains to see that $P$ acts invertibly on $C_0^\infty (L)\otimes N$. 
 If $P=t^n+r_{n-1}t^{n-1} +\cdots + r_0$, then the inverse operator is given  by
 the Taylor  expansion of $(1+x)^{-1}$, $x=r_{n-1}g^{-1} +\cdots + r_0g^{-n}$. 
\qed

We now define the {\em Grothendieck-Serre duality} functor
$\SS:D^b(Sm)^{op}\to D^b(Sm)$, $\SS:M\mapsto RHom_Z(M,D)$, where $D$ stands for the Grothendieck-Serre dualizing complex.

Let $Sm_{fg}(G)\subset Sm(G)$ be the full subcategory of finitely generated modules. 

\begin{Thm} 

For $M\in D^b(Sm_{fg})$ 
we have a canonical isomorphism $$DL\circ \SS (M)\cong RHom(M,C_c^\infty(G))[r],$$
where $r$ is the split rank of the center of $G$.
\end{Thm}

\proof  It is easy to reduce the statement to the case of a group with compact center, 
so we assume without loss of generality that $r=0$.

For every  $M\in Sm_{fg}(G)$ and an open compact
subgroup $K\subset G$ the space $M^K$ is a finitely generated module  over the center.
Thus we have a functor $i:Sm_{fg}(G)\to Sm_Z(G)$ which sends an object $M$ in $Sm$ to the same $G$-module $M$ equipped with the natural $Z$-action; it lands in $Z$-admissible modules. We use the same notation for its extension to derived categories, $i:D^b(Sm_{fg})\to D^b(Sm_Z(G))_{adm}$.

The Theorem now follows from Proposition \ref{RHom_N} and the isomorphism 
\begin{equation}\label{i!}
i^!(H\otimes D)\cong H,
\end{equation} where $i^!$ is the right adjoint to $i$.
To establish
 \eqref{i!} we use the standard isomorphism $\delta^!(\F\boxtimes D) \cong \F$ where $\delta$
is the diagonal embedding for $Spec(Z)$ and $\F\in D^b(Coh(Spec(Z))$. The isomorphism \eqref{i!} 
follows from this in view of the isomorphism $$Forg \circ i^! i \cong \delta^! \delta_* \circ Forg,$$
where $Forg:Sm\cong H-mod\to Z-mod$ is the forgetful functor, and the fact that the latter isomorphism is compatible with the action of $H$ on the functor $Forg$. 

\subsection{Serre functor}\label{Serre_f}
Recall the notion of a {\em Serre functor} on a $k$-linear triangulated category with finite dimensional 
(graded) Hom spaces \cite[Definition 3.1]{BoKa}.

Recall that for an algebra $A$ over a base field $k$ the functor $M\mapsto RHom_{A^{op}}(M^*, A^{op})$ is an endo-functor of the derived category of $A$-modules known as the {\em Nakayama functor}. 

Suppose that $A$ has finite homological dimension, and also is finite as a module
 over its center which is of finite type over $k$. Let $D^b_{fd}(A-mod)$ be the full subcategory
 in the bounded derived category consisting of modules with finite dimensional cohomology.
It is well 
known\footnote{See e.g. \cite[Example 3.2(3)]{BoKa} for the case of a finite dimensional algebra; the general case
is similar.} that in this case the Nakayama functor restricted to $D^b_{fd}(A-mod)$ is in fact a Serre functor for that category. Thus we arrive at the following
\begin{Cor} 
The functor $DL[-r]$ for $r$ as above is a Serre functor on the category $D^b(Sm(G))_{adm}$.
\end{Cor}

\begin{Rem}
A similar argument shows that the functor $DL[-r]$ on $D^b(Sm(G))$ is a relative Serre functor
over the Bernstein center, see e.g. \cite[Definition 2.5]{BeKa}.
\end{Rem}


\begin{thebibliography}{1}

\bibitem{Al} D. Alvis, {\em The duality operation in the character ring of a finite Chevalley group,} American Mathematical Society. Bulletin. New Series, {\bf 1} (1979), no 6, 907--911.

\bibitem{Aub} A.-M. Aubert, {\em Dualit\' e dans le groupe de Grothendieck de la cat\' egorie des 
repr\' esentations lisses de longueur finie dÕun groupe r\' eductif $p$-adique,}
 Trans. Amer. Math. Soc. {\bf 347} (1995), no. 6, 2179--2189,  erratum in:
  Trans. Amer. Math. Soc. {\bf 348} (1996), no. 11, 4687--4690.
 
\bibitem{cen} J. Bernstein, 
{\em Le "centre'' de Bernstein,} Edited by P. Deligne. Travaux en Cours, Representations of reductive groups over a local field, 1--32, Hermann, Paris, 1984. 


\bibitem{BBG} J. Bernstein, A. Braverman, D. Gaitsgory,
{\em The Cohen-Macaulay property of the category of $({\mathfrak g},K)$-modules,}
Selecta Math. (N.S.) {\bf 3} (1997), no. 3, 303--314. 

\bibitem{Ber2nd} J. Bernstein, {\em Second adjointness for representations of $p$-adic groups,}
preprint, available at:

\noindent
http://www.math.tau.ac.il/$\sim$bernstei/Unpublished$\underline{\  }$texts/Unpublished$\underline{\ }$list.html

\bibitem{Bez}  R. Bezrukavnikov,  {\em Homological properties of representations of $p$-adic groups related to geometry of the group at infinity,} Ph.D. thesis, arxiv preprint
 arXiv:math/0406223.
 
\bibitem{BK} R. Bezrukavnikov, D. Kazhdan, 
{\em Geometry of second adjointness for $p$-adic groups,}
 With an appendix by Y.~Varshavsky, Bezrukavnikov and Kazhdan. Represent. Theory {\bf 19} (2015), 299--332.

\bibitem{BeKa} R. Bezrukavnikov, D. Kaledin, 
\emph{McKay equivalence for symplectic resolutions of singularities,}
 Tr. Mat. Inst. Steklova  
{\bf 246}  (2004), 
 Algebr. Geom. Metody, Svyazi i Prilozh., 20--42; 
 translation in  Proc. Steklov Inst. Math.  2004,  no. 3 (246), 13--33. 

\bibitem{BoKa}  A.I. Bondal, M.M.  Kapranov, {\em Representable functors, Serre functors, and mutations,}  Math. USSR-Izv. {\bf 35} (1990), no. 3, 519--541;
translated from Izv. Akad. Nauk SSSR Ser. Mat. {\bf 53} (1989), no. 6, 1183--1205.

\bibitem{Cur} C. Curtis,
{\em Truncation and duality in the character ring of a finite group of Lie type,}
 Journal of Algebra {\bf 62} (1980) no 2, 320--332.

\bibitem{DL} P. Deligne, G. Lusztig,
{\em Duality for representations of a reductive group over a finite field,}
J. Algebra {\bf 74} (1982), no. 1, 284--291. 

\bibitem{DW} V. Drinfeld, J. Wang, 
 {\em On a strange invariant bilinear form on the space of automorphic forms,}
  Selecta Math. (N.S.) {\bf 22} (2016), no. 4, 1825--1880.
 
\bibitem{Ga} D. Gaitsgory, {\em
A "strange" functional equation for Eisenstein series and miraculous duality on the moduli stack of bundles,} preprint arXiv:1404.6780.
 
 \bibitem{Gr} M. Gromov, {\em Hyperbolic groups,} in Gersten, Steve M. ``Essays in group theory," Mathematical Sciences Research Institute Publications. 8. New York: Springer. pp. 75--263.
 
\bibitem{Hi} K. Hiraga, 
{\em On functoriality of Zelevinski involutions,}  Compos. Math.
{\bf 140} (2004), 1625--1656.

\bibitem{Ka} S. Kato,
{\em Duality for representations of a Hecke algebra,} Proc. Amer. Math. Soc.
{\bf 119} (1993), 941--946.


\bibitem{MR} I. Mirkovi\' c, S. Riche,  {\em Iwahori-Matsumoto involution and linear Koszul duality,}
 Int. Math. Res. Not. IMRN (2015), no. 1, 150--196.


\bibitem{SS}  P. Schneider, U. Stuhler, 
{\em Representation theory and sheaves on the Bruhat-Tits building,} Inst. Hautes \' Etudes Sci. Publ. Math. {\bf  85} (1997), 97--191. 

\bibitem{Ze} A. Zelevinsky,
{\em Induced representations of reductive p-adic groups, II:
On irreducible representations of $GL(n)$}, Ann. Sci.
Ecole Norm. Sup. {\bf 13}
(1980), 165--210.

 \end{thebibliography}
\end{document}